# A Stable Space-Time Finite Element Method for Parabolic Evolution Problems

**Stephen Edward Moore**

*In loving memory of Professor Francis Allotey*

**Abstract** This paper is concerned with the analysis of a new stable space-time finite element method (FEM) for the numerical solution of parabolic evolution problems in moving spatial computational domains. The discrete bilinear form is elliptic on the FEM space with respect to a discrete energy norm. This property together with a corresponding boundedness property, consistency and approximation results for the FEM spaces yield an *a priori* discretization error estimate with respect to the discrete norm. Finally, we confirm the theoretical results with numerical experiments in spatial moving domains.

**Keywords** finite element method · space-time · parabolic evolution problem · moving spatial computational domains · *a priori* discretization error estimates

**Mathematics Subject Classification (2000)** 65M12 · 65M60

## 1 Introduction

Parabolic initial-boundary value problems arise in many practical applications. For instance, heat conduction and diffusion processes but also evolution processes in life and social sciences can be modeled by parabolic evolution problems with a more general elliptic part, i.e, with convection-diffusion-reaction terms or even non-linear terms. The standard discretization methods in time and space are based on time-stepping methods combined with some spatial discretization technique like the finite element methods (FEM). The vertical method of lines discretizes first in space and then in time, whereas in the horizontal method of lines, the discretization starts with respect to the time

Katholische Hochschulgemeinde der Diözese Linz,
Petrinumstrasse 12/8/D220, A-4040, Linz.
Tel.: +43-660-6064199,
E-mail: moorekwesi@gmail.com

variable followed by the space variable, see e.g. [19]. Also time discontinuous Galerkin methods have been presented and analyzed for linear and non-linear parabolic initial boundary value problems [7]. We mention also space-time wavelet methods [15] and other space-time schemes including the $p$ and $hp$ in time versions of $hp-$finite element method to parabolic problems, see e.g., [1] and [2] respectively. This was followed by $hp$ dG time-stepping combined with FEM in space, see [14].

Another approach to discretization of parabolic evolution problems are the space-time methods. The approach is known to be a natural way of numerical discretization for problems requiring deforming and moving meshes. The space-time schemes allow for discretization in time and space simultaneously and are also applicable even for unstructured meshes see, e.g., [6]. In this approach, the time variable $t$ is considered as just another variable $x_{d+1}$ if $x_1, \ldots, x_d$ are the spatial variables. In this sense, the time derivative acts as a strong convection term in the direction $x_{d+1}$, see e.g. [17,18,3,13]. Recently, conforming space-time FEM approximations have been studied and presented in [16]. By means of an upwind-stabilized single-patch space-time Isogeometric analysis (IgA), a stable space-time scheme was developed in [10] for fixed and moving spatial computational domain. For $C^{p-1}-$ basis functions, we derived optimal approximation estimates for fixed and moving spatial computational domains with $p \geq 1$ and $p \geq 2$ respectively. For multi-patch deforming computational domains, we further derived a space-time multi-patch where dG methodology was employed in space and in time, see [12] and also [11]. The space-time approach allows for *a posteriori* error estimates, see e.g., [9].

In this paper, we generalize the results of moving spatial computational domains from [10, Section 5] to the continuous basis functions i.e. $C^0$, case. By using a time upwind test function, we derive a discrete bilinear form $a_h(\cdot, \cdot) : V_{0h} \times V_{0h} \to \mathbb{R}$ that is elliptic on $V_{0h}$ with respect to a discrete norm $\|\cdot\|_h$. The boundedness of the bilinear form is asserted in a product space $V_{h,*} \times V_h$, where the vector space $V_{h,*}$ has an associated norm $\|\cdot\|_{h,*}$. Using the ellipticity and boundedness results together with consistency results, we derive *a priori* discretization error estimates in the discrete norm $\|\cdot\|_h$. Finally, we present numerical experiments that demonstrate the feasibility of proposed scheme.

The outline of the article is as follows: In Section 2, we present the standard space-time variational scheme for the parabolic problem. Then, in Section 3, we present the new stable space-time finite element scheme. Section 4 is devoted to the derivation of *a priori* discretization error estimates. In Section 5, we present and discuss numerical experiments. We conclude the article with remarks on the scheme.

## 2 Space-time variational formulation

Let $Q \subset \mathbb{R}^{d+1}, d \in \{1, 2, 3\}$ be a bounded and Lipschitz space-time domain. Let $\alpha = (\alpha_1, ..., \alpha_d)$ be a multi-index with non-negative integers $\alpha_1, ..., \alpha_d$ and $|\alpha| = \alpha_1 + \ldots + \alpha_d$. For the multi-index $\alpha$, let $\partial_x^\alpha u := \partial^{|\alpha|} u / \partial x^\alpha =$ $\partial^{|\alpha|} u / \partial x_1^{\alpha_1} \ldots \partial x_d^{\alpha_d}$ denote the spatial derivatives with $\partial_t^i u := \partial^i u / \partial t^i$ as time derivatives. We denote the space of square-integrable functions by $L_2(Q)$ and a corresponding norm $\|\cdot\|_{L_2(Q)}$. For positive integers $s_x, s_t$, we define the Sobolev space $H^{s_x,s_t}(Q) = \{u \in L_2(Q) : \partial_x^\alpha u \in L_2(Q), \forall \alpha \text{ with } 0 \leq |\alpha| \leq s_x, \partial_t^i u \in L_2(Q), i = 0, \ldots, s_t\}$, see e.g., [8].

We consider a linear parabolic initial-boundary value problem: find $u : \overline{Q} \to \mathbb{R}$ such that

$$\partial_t u - \Delta u = f \quad \text{in} \quad Q, \quad u = 0 \quad \text{on} \quad \Sigma \cup \Sigma_0, \tag{1}$$

posed in the space-time domain $\overline{Q} := \overline{\Omega}(t) \times [0, T]$, where $\Delta$ is the Laplace operator, $\partial_t$ denotes the partial time derivative, $f$ is a given source function, $T$ is the final time, $Q := \Omega(t) \times (0, T)$, where $\Omega(t) \subset \mathbb{R}^d$ for $d = 1, 2, 3$, denotes the deforming domain $\Sigma = \partial \Omega(t) \times (0, T)$, $\Sigma_0 = \Omega(0) \times \{0\}$, $\Sigma_T = \Omega(T) \times \{T\}$ with the boundary $\partial Q := \Sigma \cup \Sigma_0 \cup \Sigma_T$.

The space-time variational formulation of (1) reads: find $u \in H^{1,0}_{0,\underline{0}}(Q)$ such that

$$a(u, v) = \ell(v) \quad \forall v \in H^{1,1}_{0,\overline{0}}(Q), \tag{2}$$

with the bilinear and linear forms given by

$$a(u, v) = -\int_Q u \partial_t v \, dx dt + \int_Q \nabla_x u \cdot \nabla_x v \, dx dt \quad \text{and} \quad \ell(v) = \int_Q f v \, dx dt, \tag{3}$$

where the trial and test spaces are defined by $H^{1,0}_{0,\underline{0}}(Q) = \{u \in L_2(Q) : \nabla_x u \in [L_2(Q)]^d, u = 0 \text{ on } \Sigma \text{ and } u = 0 \text{ on } \Sigma_0\}$ and $H^{1,1}_{0,\overline{0}}(Q) = \{u \in L_2(Q) : \nabla_x u \in [L_2(Q)]^d, \partial_t u \in L_2(Q), u = 0 \text{ on } \Sigma, \text{ and } u = 0 \text{ on } \Sigma_T\}$. We denote the gradient with respect to the spatial variables by $\nabla_x u = (\partial u / \partial x_1, \ldots, \partial u / \partial x_d)^\top$. The variational problem (2) is known to have a unique weak solution, see e.g., [8].

## 3 Stable space-time finite element method

Let $\mathcal{K}_h$ be a decomposition of the space-time domain $Q \subset \mathbb{R}^{d+1}, d = 1, 2, 3$ into non-degenerate $(d+1)$−simplices i.e. triangles if $d = 1$ or tetrahedra if $d = 2$. Let the diameter of the simplex be denoted by $h_K$ and let $h$ be defined by $h = \max\{h_K : K \in \mathcal{K}_h\}$. It is assumed that the family of triangulations is quasi-uniform, which means that there exists a positive constant $C_u$, independent of $h$ such that for all triangulations $\mathcal{K}_h$, every simplex $K \in \mathcal{K}_h$ satisfies

$$h_K \leq h \leq C_u h_K, \quad \text{for all } K \in \mathcal{K}_h. \tag{4}$$

Let $\mathbb{P}_p$ be the set of all polynomials of degree less or equal to $p$, i.e. degree $\leq p$. Then $V_h$ is the set of all continuous and piecewise polynomial functions on $Q$, i.e.

$$V_h := \left\{ v_h \in C^0(\overline{Q}) : v_h|_K \in \mathbb{P}_p(K), \quad \forall K \in \mathcal{K}_h \right\}. \tag{5}$$

Finally, we define discrete space-time space as follows

$$V_{0h} := V_h \cap H^{1,1}_{0,\underline{0}}(Q). \tag{6}$$

For every two neighboring triangles or tetrahedra $K_i, K_j \in \mathcal{K}_h$, the interior facet $F_{ij}$ is given by $F_{ij} := \overline{K}_i \cap \overline{K}_j$, $i \neq j$, if the set forms a $(d+1)$−dimensional manifold, see Figure 1. We denote the set of all interior facets of the decomposition $\mathcal{K}_h$ by $\mathcal{F}_I$ that is $\mathcal{F}_I := \left(\bigcup_{i=1}^N \partial K_i\right) \setminus \partial Q$.

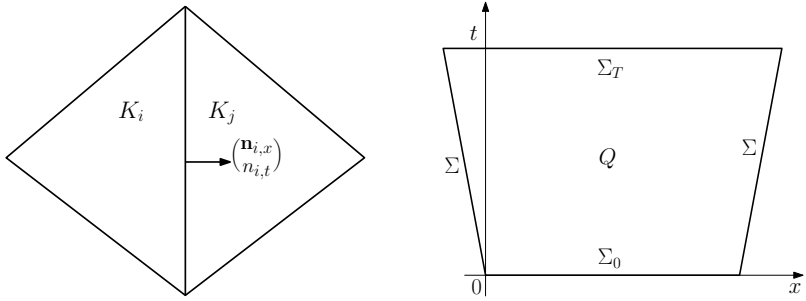

**Fig. 1** A representation of the interior facet $F_{ij}$ (left) and a description of space-time domain (right).

**Definition 1 (Jump, average, upwind)** Let $F_{ij} \in \mathcal{F}_I$ be an interior facet with the outer unit normal vector $\mathbf{n}_i = (\mathbf{n}_{i,x}, n_{i,t})^\top$ with respect to $K_i$. For a given, sufficiently smooth scalar and vector-valued functions $v$ and $\mathbf{q}$, we will denote by $v_i, \mathbf{q}_i$ the trace of the function $v$ and $\mathbf{q}$ on $F_{ij}$, and the jump across the interior facets $\mathcal{F}_I$ is defined by

$$[\![v]\!] := v_i \mathbf{n}_i + v_j \mathbf{n}_j, \quad [\![\mathbf{q}]\!] := \mathbf{q}_i \cdot \mathbf{n}_i + \mathbf{q}_j \cdot \mathbf{n}_j,$$

The jump in space direction is given by

$$[\![v]\!]_x := v_i \mathbf{n}_{i,x} + v_j \mathbf{n}_{j,x}, \quad [\![\mathbf{q}]\!]_x := \mathbf{q}_i \cdot \mathbf{n}_{i,x} + \mathbf{q}_j \cdot \mathbf{n}_{j,x},$$

whereas the jump in time direction is defined by

$$[\![v]\!]_t := v_i n_{i,t} + v_j n_{j,t}, \qquad [\![\mathbf{q}]\!]_t := \mathbf{q}_i n_{i,t} + \mathbf{q}_j n_{j,t},$$

The average of a function $v$ on the interior facet $F_{ij}$ is nothing but

$$\{v\} := \frac{1}{2}(v_i + v_j), \quad \{\mathbf{q}\} := \frac{1}{2}(\mathbf{q}_i + \mathbf{q}_j),$$

and the upwind in time direction is given by

$$\{v\}^{up} := \begin{cases} v_i & \text{for } n_{i,t} \geq 0, \\ v_j & \text{for } n_{i,t} < 0, \end{cases} \quad \{\mathbf{q}\}^{up} := \begin{cases} \mathbf{q}_i & \text{for } n_{i,t} \geq 0, \\ \mathbf{q}_j & \text{for } n_{i,t} < 0, \end{cases} \tag{7}$$

whereas the downwind in time direction is given by

$$\{v\}^{down} := \begin{cases} v_j & \text{for } n_{i,t} \geq 0, \\ v_i & \text{for } n_{i,t} < 0, \end{cases} \quad \{\mathbf{q}\}^{down} := \begin{cases} \mathbf{q}_j & \text{for } n_{i,t} \geq 0, \\ \mathbf{q}_i & \text{for } n_{i,t} < 0. \end{cases} \tag{8}$$

The next lemma provides us with some elementary formulas that are necessary for the derivation of our variational space-time FEM scheme.

**Lemma 1** *Let $F_{ij} \in \mathcal{F}_I$ be an interface, and let $u$ and $v$ be sufficiently smooth functions on the interface. Then the following formulas hold:*

$$[\![uv]\!]_x = \{u\}[\![v]\!]_x + \{v\}[\![u]\!]_x, \tag{9}$$

$$[\![uv]\!]_t = \{u\}^{up}[\![v]\!]_t + \{v\}^{down}[\![u]\!]_t. \tag{10}$$

*Proof* See e.g. [11]. □

The derivation of the scheme is as follows: we multiply our model problem (1) with a test function of the form $v_h + \theta h \partial_t v_h$ for an arbitrary $v_h \in V_{0h} \subset H_{0,\underline{0}}^{1,1}(Q)$ and a positive constant $\theta$ which will be determined later,

$$\int_Q f(v_h + \theta h \partial_t v_h)\,dxdt = \int_Q \left(\partial_t u(v_h + \theta h \partial_t v_h) - \Delta u v_h - \theta h \Delta u \partial_t v_h\right) dxdt. \tag{11}$$

Integration by parts with respect to $x$ in the second term of the bilinear form on the right-hand side of (11) gives

$$\int_Q f(v_h + \theta h \partial_t v_h)\,dxdt$$
$$= \int_Q \left(\partial_t u(v_h + \theta h \partial_t v_h) + \nabla_x u \cdot \nabla_x v_h - \theta h \Delta u \partial_t v_h\right) dxdt - \int_{\partial Q} \mathbf{n}_x \cdot \nabla_x u\, v_h\, ds,$$

where using the facts that $v_h \in V_{0h}$ and $\mathbf{n}_x$ are zero on $\Sigma_0$ and $\Sigma_T$ yields

$$\int_Q f(v_h + \theta h \partial_t v_h)\,dxdt$$
$$= \int_Q \left(\partial_t u(v_h + \theta h \partial_t v_h) + \nabla_x u \cdot \nabla_x v_h - \theta h \Delta u \partial_t v_h\right) dxdt.$$

Concerning the last term, we sum over each element and apply an integration by parts with respect to the spatial direction

$$-\int_Q \Delta u \partial_t v_h\,dxdt = -\sum_{K \in \mathcal{K}_h} \int_K \Delta u \partial_t v_h\,dxdt$$
$$= \sum_{K \in \mathcal{K}_h} \int_K \nabla_x u \cdot \nabla_x \partial_t v_h\,dxdt - \sum_{K \in \mathcal{K}_h} \int_{\partial K} \mathbf{n}_x \cdot \nabla_x u \partial_t v_h\,ds$$
$$= \sum_{K \in \mathcal{K}_h} \int_K \nabla_x u \cdot \nabla_x \partial_t v_h\,dxdt - \sum_{F_{ij} \in \mathcal{F}_I} \int_{F_{ij}} [\![\nabla_x u \partial_t v_h]\!]_x\,ds - \int_{\Sigma} \mathbf{n}_x \cdot \nabla_x u \partial_t v_h\,ds.$$

By an integration by parts with respect to the time direction for the first term above and using the fact that $v_h \in V_{0h}$, we obtain

$$-\int_Q \Delta u \partial_t v_h\,dxdt = -\sum_{K \in \mathcal{K}_h} \int_K \partial_t \nabla_x u \cdot \nabla_x v_h\,dxdt + \sum_{F_{ij} \in \mathcal{F}_I} \int_{F_{ij}} [\![\nabla_x u \cdot \nabla_x v_h]\!]_t\,ds$$

$$+ \int_{\Sigma \cup \Sigma_T} n_t \nabla_x u \cdot \nabla_x v_h \, ds - \sum_{F_{ij} \in \mathcal{F}_I} \int_{F_{ij}} [\![\nabla_x u \partial_t v_h]\!]_x \, ds - \int_{\Sigma} \mathbf{n}_x \cdot \nabla_x u \partial_t v_h \, ds.$$

Now, considering the terms on the interior facets $F_{ij} \in \mathcal{F}_I$, we use the identities (9) and (10) from Lemma 1 to obtain

$$\sum_{F_{ij} \in \mathcal{F}_I} \int_{F_{ij}} \left( [\![\nabla_x u \cdot \nabla_x v_h]\!]_t - [\![\nabla_x u \partial_t v_h]\!]_x \right) ds$$

$$= \sum_{F_{ij} \in \mathcal{F}_I} \int_{F_{ij}} \left( \{\nabla_x u\}^{up} \cdot [\![\nabla_x v_h]\!]_t + \{\nabla_x v_h\}^{down} \cdot [\![\nabla_x u]\!]_t \right.$$

$$\left. - \{\nabla_x u\} \cdot [\![\partial_t v_h]\!]_x - \{\partial_t v_h\} [\![\nabla_x u]\!]_x \right) ds.$$

Assuming that our solution $u$ belongs to $H^2(Q)$, the jumps of the derivative of $u$ are zero i.e. $[\![\nabla_x u]\!]_x = 0$ and $[\![\nabla_x u]\!]_t = 0$, thus yielding

$$\sum_{F_{ij} \in \mathcal{F}_I} \int_{F_{ij}} [\![\nabla_x u \cdot \nabla_x v_h]\!]_t - [\![\nabla_x u \partial_t v_h]\!]_x \, ds$$

$$= \sum_{F_{ij} \in \mathcal{F}_I} \int_{F_{ij}} (\{\nabla_x u\}^{up} \cdot [\![\nabla_x v_h]\!]_t - \{\nabla_x u\} \cdot [\![\partial_t v_h]\!]_x) ds.$$

Also, considering the terms on the boundary $\Sigma$, we have

$$\int_{\Sigma} \nabla_x u \cdot \left( n_t \nabla_x v_h - \mathbf{n}_x \partial_t v_h \right) ds = 0,$$

since $(n_t \nabla_x v_h - \mathbf{n}_x \partial_t v_h)$ are the tangential derivatives of $v_h$ and $v_h = 0$ on $\Sigma$ as discussed in [10]. Since the solution $u$ is assumed to be smooth enough, we have that $[\![\partial_t u]\!]_x = 0$. This enables us to add the consistent term

$$\theta h \sum_{F_{ij} \in \mathcal{F}_I} \int_{F_{ij}} \left( \{\nabla_x v_h\} \cdot [\![\partial_t u]\!]_x + \delta [\![\partial_t u]\!]_x \cdot [\![\partial_t v_h]\!]_x \right) ds, \qquad (12)$$

to the bilinear form where $\delta$ is a positive constant. Finally, we can write the variational space-time FEM scheme as follows: find $u_h \in V_{0h}$ such that

$$a_h(u_h, v_h) = \ell_h(v_h) \quad \forall v_h \in V_{0h}, \qquad (13)$$

where the bilinear form is given by

$$a_h(u_h, v_h) = \int_Q \left( \partial_t u_h (v_h + \theta h \partial_t v_h) + \nabla_x u \cdot \nabla_x v_h \right) dx dt$$

$$- \theta h \sum_{K \in \mathcal{K}_h} \int_K \partial_t \nabla_x u_h \cdot \nabla_x v_h \, dx dt + \theta h \int_{\Sigma_T} \nabla_x u_h \cdot \nabla_x v_h \, ds$$

$$+ \theta h \sum_{F_{ij} \in \mathcal{F}_I} \int_{F_{ij}} \left( \{\nabla_x u_h\}^{up} \cdot [\![\nabla_x v_h]\!]_t - \{\nabla_x u_h\} \cdot [\![\partial_t v_h]\!]_x + \{\nabla_x v_h\} \cdot [\![\partial_t u_h]\!]_x \right.$$

$$+ \delta [\![\partial_t u_h]\!]_x \cdot [\![\partial_t v_h]\!]_x \bigg) \, ds \tag{14}$$

and the linear form

$$\ell_h(v_h) = \int_Q f\,(v_h + \theta h \partial_t v_h)\,dxdt. \tag{15}$$

*Remark 1* Although, we have considered homogeneous Dirichlet boundary and initial conditions, it is also possible to consider other boundary conditions including Neumann and Robin boundary conditions as well as non-homogeneous Dirichlet boundary conditions. Particularly, the weak imposition of the Dirichlet and initial boundary conditions will be an interesting consideration.

**Theorem 1** *If the solution $u \in H^{1,0}_{0,\underline{0}}(Q)$ of the variational problem (2) belongs to $H^2(Q)$ then the solution satisfies the consistency identity*

$$a_h(u, v_h) = \ell_h(v_h), \quad \forall\, v_h \in V_{0h}. \tag{16}$$

An immediate consequence of consistency is the Galerkin orthogonality property

$$a_h(u - u_h, v_h) = 0, \quad \forall v_h \in V_{0h}. \tag{17}$$

Due to the assumptions on the mesh elements, the following Lemmata that are required for the analysis of the scheme hold. We refer the interested reader to [4] for the proofs.

**Lemma 2** *Let $K \in \mathcal{K}_h$. Then the scaled trace inequality*

$$\|v\|_{L_2(\partial K)} \le C_{tr} h_K^{-1/2} \left( \|v\|_{L_2(K)} + h_K |v|_{H^1(K)} \right), \tag{18}$$

*holds for all $v \in H^1(K)$, where $h_K$ denotes the mesh size of $K$ and $C_{tr}$ is a positive constant.*

**Lemma 3** *Let $K$ be an arbitrary mesh element from $\mathcal{K}_h$. Then the inverse inequalities*

$$\|\nabla_x v_h\|_{L_2(K)} \le C_{inv,1} h_K^{-1} \|v_h\|_{L_2(K)}, \tag{19}$$

$$\|v_h\|_{L_2(\partial K)} \le C_{inv,0} h_K^{-1/2} \|v_h\|_{L_2(K)}, \tag{20}$$

*hold for all $v_h \in V_h$, where $C_{inv,1}$ and $C_{inv,0}$ are positive constants.*

To show the coercivity of the bilinear form, we will need the following lemma.

**Lemma 4** *Let $F_{ij} \in \mathcal{F}_I$ be an interior facet and $v \in V_{0h}$ be a function. Then the following identity holds:*

$$\{v\}^{up}[\![v]\!]_t - \frac{1}{2}[\![v^2]\!]_t = \frac{1}{2}|n_{i,t}|[\![v]\!]^2. \tag{21}$$

*Proof* See e.g. [11] or [12]. □

Next, we show that the discrete bilinear form (14) is $V_{0h}-$coercive with respect to the mesh-dependent norm

$$\|v_h\|_h := \bigg(\|\nabla_x v_h\|^2_{L_2(Q)} + \theta h\|\partial_t v_h\|^2_{L_2(Q)} + \frac{1}{2}\|v_h\|^2_{L_2(\Sigma_T)} + \frac{\theta h}{2}\|\nabla_x v_h\|^2_{L_2(\Sigma_T)}$$
$$+ \frac{\theta h}{2}\sum_{F_{ij} \in \mathcal{F}_I} \|[\![\nabla_x v_h]\!]_t\|^2_{L_2(F_{ij})} + \delta\theta h \sum_{F_{ij} \in \mathcal{F}_I} \|[\![\partial_t v_h]\!]_x\|^2_{L_2(F_{ij})}\bigg)^{1/2}. \quad (22)$$

*Remark 2* The equation above (22) is a norm. Indeed, if $\|v_h\|_h = 0$ for some $v_h \in V_{0h}$, then $\nabla_x v_h = 0$ and $\partial_t v_h = 0$ in $Q$, i.e., $v_h$ is a constant in $Q$. Furthermore, $v_h \in V_{0h}$ implies that $v_h$ is zero on $\Sigma$ and $\Sigma_0$, i.e. this constant must be zero. Therefore, $v_h = 0$ in the whole space-time computational domain $Q$. The other norm axioms are trivial.

**Lemma 5** *Let $\theta > 0$ be sufficiently small such that $\theta \leq (C^2_{inv,0}C^2_u)^{-1}$, where $C_{inv,0}$ is chosen as in Lemma 3 and $C_u$ as in Assumption (4). Then the discrete bilinear form $a_h(\cdot,\cdot) : V_{0h} \times V_{0h} \to \mathbb{R}$, defined by (14), is $V_{0h}-$coercive with respect to the norm $\|\cdot\|_h$, i.e.*

$$a_h(v_h, v_h) \geq \mu_c \|v_h\|^2_h, \quad \forall v_h \in V_{0h}, \quad (23)$$

*with $\mu_c = 1/2$.*

*Proof* Let $v_h = u_h$ in (14), then we get

$$a_h(v_h, v_h) = \int_Q \big(\partial_t v_h v_h + \theta h(\partial_t v_h)^2 + |\nabla_x v_h|^2\big)\, dxdt + \theta h\int_{\Sigma_T} |\nabla_x v_h|^2\, ds$$
$$- \theta h \sum_{K \in \mathcal{K}_h} \int_K \partial_t \nabla_x v_h \cdot \nabla_x v_h\, dxdt + \theta h \sum_{F_{ij} \in \mathcal{F}_I} \int_{F_{ij}} \{\nabla_x v_h\}^{up} \cdot [\![\nabla_x v_h]\!]_t\, ds$$
$$- \theta h \sum_{F_{ij} \in \mathcal{F}_I} \int_{F_{ij}} \{\nabla_x v_h\} \cdot [\![\partial_t v_h]\!]_x\, ds + \theta h \sum_{F_{ij} \in \mathcal{F}_I} \int_{F_{ij}} \{\nabla_x v_h\} \cdot [\![\partial_t v_h]\!]_x\, ds$$
$$+ \delta\theta h \sum_{F_{ij} \in \mathcal{F}_I} \int_{F_{ij}} [\![\partial_t v_h]\!]^2_x\, ds.$$

Using Gauss' theorem, we obtain

$$a_h(v_h, v_h) = \frac{1}{2}\int_Q \partial_t v_h^2 + \theta h\|\partial_t v_h\|^2_{L_2(Q)} + \|\nabla_x v_h\|^2_{L_2(Q)} + \theta h\|\nabla_x v_h\|^2_{L_2(\Sigma_T)}$$
$$- \frac{\theta h}{2} \sum_{K \in \mathcal{K}_h} \int_K \partial_t |\nabla_x v_h|^2\, dxdt + \theta h \sum_{F_{ij} \in \mathcal{F}_I} \int_{F_{ij}} \{\nabla_x v_h\}^{up} \cdot [\![\nabla_x v_h]\!]_t\, ds$$
$$+ \delta\theta h \sum_{F_{ij} \in \mathcal{F}_I} \int_{F_{ij}} [\![\partial_t v_h]\!]^2\, ds$$

$$= \frac{1}{2} \int_{\partial Q} n_t v_h^2 + \theta h \|\partial_t v_h\|_{L_2(Q)}^2 + \|\nabla_x v_h\|_{L_2(Q)}^2 + \theta h \|\nabla_x v_h\|_{L_2(\Sigma_T)}^2$$
$$- \frac{\theta h}{2} \sum_{K \in \mathcal{K}_h} \int_K \partial_t |\nabla_x v_h|^2 \, dxdt + \theta h \sum_{F_{ij} \in \mathcal{F}_I} \int_{F_{ij}} \{\nabla_x v_h\}^{up} \cdot [\![\nabla_x v_h]\!]_t \, ds$$
$$+ \delta\theta h \sum_{F_{ij} \in \mathcal{F}_I} \|[\![\partial_t v_h]\!]_x\|_{L_2(F_{ij})}^2.$$

Rewriting the boundary terms such that $\partial Q := \Sigma \cup \Sigma_0 \cup \Sigma_T$ and $\mathcal{F}_I := \left(\bigcup_{i=1}^N \partial K_i\right) \setminus \partial Q$ with the interior facet $F_{ij} \subset \mathcal{F}_I$ and using $v_h = 0$ on $\Sigma_0$ yields

$$a_h(v_h, v_h) = \frac{1}{2} \int_{\Sigma_T} v_h^2 \, ds + \theta h \|\partial_t v_h\|_{L_2(Q)}^2 + \|\nabla_x v_h\|_{L_2(Q)}^2 + \theta h \|\nabla_x v_h\|_{L_2(\Sigma_T)}^2$$
$$- \frac{\theta h}{2} \sum_{F_{ij} \in \mathcal{F}_I} \int_{F_{ij}} [\![|\nabla_x v_h|^2]\!]_t \, ds - \frac{\theta h}{2} \|\nabla_x v_h\|_{L_2(\Sigma_T)}^2$$
$$- \frac{\theta h}{2} \int_\Sigma n_t |\nabla_x v_h|^2 \, ds + \theta h \sum_{F_{ij} \in \mathcal{F}_I} \int_{F_{ij}} \{\nabla_x v_h\}^{up} \cdot [\![\nabla_x v_h]\!]_t \, ds$$
$$+ \delta\theta h \sum_{F_{ij} \in \mathcal{F}_I} \|[\![\partial_t v_h]\!]_x\|_{L_2(F_{ij})}^2.$$

Using identity (21) and the fact that $v_h \in V_{0h}$, we obtain

$$a_h(v_h, v_h) = \frac{1}{2}\|v_h\|_{L_2(\Sigma_T)}^2 + \theta h \|\partial_t v_h\|_{L_2(Q)}^2 + \|\nabla_x v_h\|_{L_2(Q)}^2 + \frac{\theta h}{2}\|\nabla_x v_h\|_{L_2(\Sigma_T)}^2$$
$$- \frac{\theta h}{2} \int_\Sigma n_t |\nabla_x v_h|^2 \, ds + \theta h \sum_{F_{ij} \in \mathcal{F}_I} \int_{F_{ij}} \left(\{\nabla_x v_h\}^{up} \cdot [\![\nabla_x v_h]\!]_t \, ds - \frac{1}{2}[\![|\nabla_x v_h|^2]\!]_t\right) ds$$
$$+ \delta\theta h \sum_{F_{ij} \in \mathcal{F}_I} \|[\![\partial_t v_h]\!]_x\|_{L_2(F_{ij})}^2$$
$$= \frac{1}{2}\|v_h\|_{L_2(\Sigma_T)}^2 + \theta h \|\partial_t v_h\|_{L_2(Q)}^2 + \|\nabla_x v_h\|_{L_2(Q)}^2 + \frac{\theta h}{2}\|\nabla_x v_h\|_{L_2(\Sigma_T)}^2$$
$$- \frac{\theta h}{2} \int_\Sigma n_t |\nabla_x v_h|^2 \, ds + \frac{\theta h}{2} \sum_{F_{ij} \in \mathcal{F}_I} \int_{F_{ij}} |n_{i,t}| \left([\![|\nabla_x v_h|]\!]\right)^2 ds + \delta\theta h \sum_{F_{ij} \in \mathcal{F}_I} \|[\![\partial_t v_h]\!]_x\|_{L_2(F_{ij})}^2.$$

By using $|n_{i,t}| \geq |n_{i,t}|^2$ and Definition 1, we obtain

$$a_h(v_h, v_h) \geq \frac{1}{2}\|v_h\|_{L_2(\Sigma_T)}^2 + \theta h \|\partial_t v_h\|_{L_2(Q)}^2 + \|\nabla_x v_h\|_{L_2(Q)}^2 + \frac{\theta h}{2}\|\nabla_x v_h\|_{L_2(\Sigma_T)}^2$$
$$- \frac{\theta h}{2}\|\nabla_x v_h\|_{L_2(\Sigma)}^2 + \frac{\theta h}{2} \sum_{F_{ij} \in \mathcal{F}_I} \|[\![\nabla_x v_h]\!]_t\|_{L_2(F_{ij})}^2 + \delta\theta h \sum_{F_{ij} \in \mathcal{F}_I} \|[\![\partial_t v_h]\!]_x\|_{L_2(F_{ij})}^2$$
$$\geq \|v_h\|_h^2 - \frac{\theta h}{2}\|\nabla_x v_h\|_{L_2(\Sigma)}^2. \qquad (24)$$

Now by applying (20) and the quasi-uniform assumption (4) to the term on $\Sigma$, we have

$$\theta h/2 \|\nabla_x v_h\|^2_{L_2(\Sigma)} \leq \theta h/2 C^2_{inv,0} C^2_u h^{-1} \|\nabla_x v_h\|^2_{L_2(Q)} = (\theta C^2_{inv,0} C^2_u)/2 \|\nabla_x v_h\|^2_{L_2(Q)}. \tag{25}$$

By inserting (25) into (24), we have

$$\begin{aligned} a_h(v_h, v_h) &\geq \|v_h\|^2_h - (\theta C^2_{inv,0} C^2_u)/2 \|\nabla_x v_h\|^2_{L_2(Q)}, \\ &\geq (1 - (\theta C^2_{inv,0} C^2_u)/2) \|v_h\|^2_h \geq \frac{1}{2}\|v_h\|^2_h, \end{aligned}$$

since $\theta \leq (C^2_{inv,0} C^2_u)^{-1}$. This completes the proof of the coercivity. □

*Remark 3* We observe that Lemma 5 is sufficient for well-posedness of the discrete problem (13). The Lemma 5 implies the uniqueness of the solution $u_h \in V_{0h}$. Since the space-time FEM scheme (13) is posed in a finite-dimensional space $V_{0h}$, the uniqueness yields the existence of the solution.

We need uniform boundedness of the discrete bilinear form $a_h(\cdot,\cdot)$ on $V_{0h,*} \times V_{0h}$, where the space $V_{0h,*} = H^{1,0}_0(Q) \cap H^2(Q) + V_{0h}$ is equipped with the norm

$$\begin{aligned} \|v\|_{h,*} := \Bigg( &\|\nabla_x v\|^2_{L_2(Q)} + \theta h \|\partial_t v\|^2_{L_2(Q)} + \frac{1}{2}\|v\|^2_{L_2(\Sigma_T)} + \theta h \|\nabla_x v\|^2_{L_2(\Sigma_T)} \\ &+ \theta h \sum_{F_{ij} \in \mathcal{F}_I} \|[\![\nabla_x v]\!]_t\|^2_{L_2(F_{ij})} + \delta \theta h \sum_{F_{ij} \in \mathcal{F}_I} \|[\![\partial_t v_h]\!]_x\|^2_{L_2(F_{ij})} \\ &+ (\theta h)^{-1} \|v\|^2_{L_2(Q)} + (\theta h)^2 \|\partial_t \nabla_x v\|^2_{L_2(Q)} + \theta h \sum_{F_{ij} \in \mathcal{F}_I} \|\{\nabla_x v\}^{up}\|^2_{L_2(F_{ij})} \\ &+ \sum_{i=1}^N h_{K_i} \|\partial_t v\|^2_{L_2(\partial K_i)} \Bigg)^{\frac{1}{2}}. \end{aligned} \tag{26}$$

**Lemma 6** *The discrete bilinear form $a_h(\cdot,\cdot)$, defined by (14), is uniformly bounded on $V_{0h,*} \times V_{0h}$, i.e., there exists a positive constant $\mu_b$ that does not depend on h such that*

$$|a_h(u, v_h)| \leq \mu_b \|u\|_{h,*} \|v_h\|_h, \ \forall\, u \in V_{0h,*}, \forall\, v_h \in V_{0h}. \tag{27}$$

*with $\mu_b = \max\{\delta^{-1} C^{-2}_{inv,0}/2, 4 + \delta^{-1}/2\}^{1/2}$, where $C_{inv,0}$ is chosen as in Lemma 3, $\delta$ is a positive constant and $\theta$ chosen as in Lemma 5.*

*Proof* Let us estimate the bilinear form (14) as follows: for the first term, we apply integration by parts and Cauchy-Schwarz inequality yielding

$$\begin{aligned} \int_Q \partial_t u v_h \, dxdt &= -\int_Q u \partial_t v_h \, dxdt + \int_{\partial Q} u n_t v_h \, ds \\ &\leq \left((\theta h)^{-1} \|u\|^2_{L_2(Q)}\right)^{1/2} \left((\theta h) \|\partial_t v_h\|^2_{L_2(Q)}\right)^{1/2} \end{aligned}$$

$$+ \left(\|u\|_{L_2(\Sigma_T)}^2\right)^{1/2} \left(\|v_h\|_{L_2(\Sigma_T)}^2\right)^{1/2}. \tag{28}$$

By using Cauchy Schwarz's inequality on the following terms, we obtain

$$\theta h \int_{\Sigma_T} \nabla_x u \cdot \nabla_x v_h \, ds \leq \left(\theta h \|\nabla_x u\|_{L_2(\Sigma_T)}^2\right)^{1/2} \left(\theta h \|\nabla_x v_h\|_{L_2(\Sigma_T)}^2\right)^{1/2},$$

$$\theta h \int_Q \partial_t u \partial_t v_h \, dxdt \leq \left(\theta h \|\partial_t u\|_{L_2(Q)}^2\right)^{1/2} \left(\theta h \|\partial_t v_h\|_{L_2(Q)}^2\right)^{1/2},$$

$$\int_Q \nabla_x u \cdot \nabla_x v_h \, dxdt \leq \left(\|\nabla_x u\|_{L_2(Q)}^2\right)^{1/2} \left(\|\nabla_x v_h\|_{L_2(Q)}^2\right)^{1/2},$$

$$\theta h \sum_{K \in \mathcal{K}_h} \int_K \partial_t \nabla_x u \cdot \nabla_x v_h \, dxdt \leq \left((\theta h)^2 \|\partial_t \nabla_x u\|_{L_2(Q)}^2\right)^{1/2} \left(\|\nabla_x v_h\|_{L_2(Q)}^2\right)^{1/2},$$

$$\delta \theta h \sum_{F_{ij} \in \mathcal{F}_I} \int_{F_{ij}} [\![\partial_t u]\!]_x \cdot [\![\partial_t v_h]\!]_x \, ds \leq \left(\delta \theta h \sum_{F_{ij} \in \mathcal{F}_I} \|[\![\partial_t u]\!]_x\|_{L_2(F_{ij})}^2\right)^{1/2}$$
$$\times \left(\delta \theta h \sum_{F_{ij} \in \mathcal{F}_I} \|[\![\partial_t v_h]\!]_x\|_{L_2(F_{ij})}^2\right)^{1/2}.$$

Next, we estimate the interface terms using the quasi-uniform assumption (4) and the inverse inequality (20)

$$\theta h \sum_{F_{ij} \in \mathcal{F}_I} \int_{F_{ij}} \{\nabla_x u\}^{up} \cdot [\![\nabla_x v_h]\!]_t \, ds$$
$$\leq \left(\theta h \sum_{F_{ij} \in \mathcal{F}_I} \|\{\nabla_x u\}^{up}\|_{L_2(F_{ij})}^2\right)^{1/2} \left(\theta h \sum_{F_{ij} \in \mathcal{F}_I} \|[\![\nabla_x v_h]\!]_t\|_{L_2(F_{ij})}^2\right)^{1/2}$$
$$\leq \left(\theta h \sum_{F_{ij} \in \mathcal{F}_I} \|\{\nabla_x u\}^{up}\|_{L_2(F_{ij})}^2\right)^{1/2} \left(2\theta C_u \sum_{i=1}^N h_{K_i} \|\nabla_x v_{h,i}\|_{L_2(\partial K_i)}^2\right)^{1/2}$$
$$\leq \left(\theta h \sum_{F_{ij} \in \mathcal{F}_I} \|\{\nabla_x u\}^{up}\|_{L_2(F_{ij})}^2\right)^{1/2} \left(2\theta C_u C_{inv,0}^2 \|\nabla_x v_h\|_{L_2(Q)}^2\right)^{1/2}, \tag{29}$$

and

$$\theta h \sum_{F_{ij} \in \mathcal{F}_I} \int_{F_{ij}} \{\nabla_x v_h\} \cdot [\![\partial_t u]\!]_x \, ds$$
$$\leq \left(\delta^{-1} \theta h \sum_{F_{ij} \in \mathcal{F}_I} \|\{\nabla_x v_h\}\|_{L_2(F_{ij})}^2\right)^{1/2} \left(\delta \theta h \sum_{F_{ij} \in \mathcal{F}_I} \|[\![\partial_t u]\!]_x\|_{L_2(F_{ij})}^2\right)^{1/2}$$
$$\leq \left(\delta^{-1} \frac{\theta}{2} C_u \sum_{i=1}^N h_{K_i} \|\nabla_x v_{h,i}\|_{L_2(\partial K_i)}^2\right)^{1/2} \left(\delta \theta h \sum_{F_{ij} \in \mathcal{F}_I} \|[\![\partial_t u]\!]_x\|_{L_2(F_{ij})}^2\right)^{1/2}$$

$$\leq \left(\delta^{-1}\frac{\theta}{2}C_u C_{inv,0}^2 \|\nabla_x v_h\|_{L_2(Q)}^2\right)^{1/2} \left(\delta\theta h \sum_{F_{ij}\in\mathcal{F}_I} \|[\![\partial_t u]\!]_x\|_{L_2(F_{ij})}^2\right)^{1/2}. \quad (30)$$

The final term is estimated using the quasi-uniform assumption (4) to arrive at

$$\theta h \sum_{F_{ij}\in\mathcal{F}_I} \int_{F_{ij}} \{\nabla_x u\} \cdot [\![\partial_t v_h]\!]_x \, ds$$

$$\leq \left(\delta^{-1}\theta h \sum_{F_{ij}\in\mathcal{F}_I} \|\{\nabla_x u\}\|_{L_2(F_{ij})}^2\right)^{1/2} \left(\delta\theta h \sum_{F_{ij}\in\mathcal{F}_I} \|[\![\partial_t v_h]\!]_x\|_{L_2(F_{ij})}^2\right)^{1/2}$$

$$\leq \left(\delta^{-1}\theta C_u/2 \sum_{i=1}^N h_{K_i} \|\nabla_x u_i\|_{L_2(\partial K_i)}^2\right)^{1/2} \left(\delta\theta h \sum_{F_{ij}\in\mathcal{F}_I} \|[\![\partial_t v_h]\!]_x\|_{L_2(F_{ij})}^2\right)^{1/2}. \quad (31)$$

Finally, substituting equations (29), (30) and (31) into (14), we get

$$|a_h(u, v_h)| \leq \Big((\theta h)^{-1}\|u\|_{L_2(Q)}^2 + \|u\|_{L_2(\Sigma_T)}^2 + \theta h\|\partial_t u\|_{L_2(Q)}^2 + \|\nabla_x u\|_{L_2(Q)}^2$$

$$+ \theta h\|\nabla_x u\|_{L_2(\Sigma_T)}^2 + (\theta h)^2\|\partial_t \nabla_x u\|_{L_2(Q)}^2 + \theta h \sum_{F_{ij}\in\mathcal{F}_I} \|\{\nabla_x u\}^{up}\|_{L_2(F_{ij})}^2$$

$$+ \delta^{-1}\theta C_u/2 \sum_{i=1}^N h_{K_i}\|\nabla_x u_i\|_{L_2(\partial K_i)}^2 + 2\delta\theta h \sum_{F_{ij}\in\mathcal{F}_I} \|[\![\partial_t u]\!]_x\|_{L_2(F_{ij})}^2\Big)^{1/2}$$

$$\times \Big(\theta h\|\partial_t v_h\|_{L_2(Q)}^2 + \|v_h\|_{L_2(\Sigma_T)}^2 + \theta h\|\partial_t v_h\|_{L_2(Q)}^2 + 2\|\nabla_x v_h\|_{L_2(Q)}^2$$

$$+ \theta h\|\nabla_x v_h\|_{L_2(\Sigma_T)}^2 + 2\theta C_u C_{inv,0}^2\|\nabla_x v_h\|_{L_2(Q)}^2$$

$$+ (\delta^{-1}\theta C_u C_{inv,0}^2)/2\|\nabla_x v_h\|_{L_2(Q)}^2 + 2\delta\theta h \sum_{F_{ij}\in\mathcal{F}_I} \|[\![\partial_t v_h]\!]_x\|_{L_2(F_{ij})}^2\Big)^{1/2}$$

$$\leq \mu_b \|u\|_{h,*} \|v_h\|_h,$$

with $\mu_b = \max\left\{2(1 + C_u C_{inv,0}^2 + \delta^{-1}\theta C_u C_{inv,0}^2/4), \delta^{-1}\theta C_u/2\right\}^{1/2}$ where $\theta \leq \left(C_u C_{inv,0}^2\right)^{-1}$. □

## 4 A priori discretization error estimates

To derive *a priori* error estimates, we will require some interpolation results for finite element. Let us recall the following interpolation result, see e.g. [4].

**Lemma 7** *Let $l, d$ and $s$ be positive integers such that $2 \leq s \leq p+1$ with $s > (d+1)/2$ and $0 \leq l \leq s$. Furthermore, let $v \in H^s(Q) \cap H_{0,\underline{0}}^{1,1}(Q)$. Then* there exists a projective operator $\Pi_h$ from $H_{0,\underline{0}}^{1,1}(Q) \cap H^s(Q)$ to $V_{0h}$ and a positive generic constant $C_s$ such that

$$\left( \sum_{K \in \mathcal{K}_h} |v - \Pi_h v|^2_{H^l(K)} \right)^{1/2} \leq C_s h^{s-l} \|v\|_{H^s(Q)}, \tag{32}$$

where h denotes the maximum mesh-size parameter in the physical domain and the constant $C_s$ only depends on $l, s, p$ and the shape regularity of the physical domain Q.

The interpolation result enables us to derive estimates for terms on the interior facets.

**Lemma 8** *Let the assumptions of Lemma 7 hold. Then the following estimates hold*

$$\sum_{F_{ij} \in \mathcal{F}_I} h \|[\![\nabla_x(v - \Pi_h v)]\!]_t\|^2_{L_2(F_{ij})} \leq C_3 h^{2(r-1)} \|v\|^2_{H^r(Q)}, \tag{33}$$

$$\sum_{F_{ij} \in \mathcal{F}_I} h \|\{\nabla_x(v - \Pi_h v)\}^{up}\|^2_{L_2(F_{ij})} \leq C_4 h^{2(r-1)} \|v\|^2_{H^r(Q)}, \tag{34}$$

$$\sum_{F_{ij} \in \mathcal{F}_I} \delta\theta h \|[\![\partial_t(v - \Pi_h v)]\!]_x\|^2_{L_2(F_{ij})} \leq C_5 h^{2(r-1)} \|v\|^2_{H^r(Q)}, \tag{35}$$

*where $r = \min\{s, p+1\}$, p denotes the underlying polynomial degree, $\delta$ is a positive constant and the generic constants $C_3, C_4$ and $C_5$ do not depend on the mesh element size h.*

*Proof* By using the quasi-uniformity assumption (4), the trace inequality (18), and the approximation property (32), we obtain

$$\sum_{F_{ij} \in \mathcal{F}_I} h \|[\![\nabla_x(v - \Pi_h v)]\!]_t\|^2_{L_2(F_{ij})}$$

$$\leq 2 C_u C_{tr}^2 \sum_{i=1}^N h_{K_i} h_{K_i}^{-1} \left( \|\nabla_x(v - \Pi_{h,i} v)\|^2_{L_2(K_i)} + h_{K_i}^2 |\nabla_x(v - \Pi_{h,i} v)|^2_{H^1(K_i)} \right)$$

$$\leq 4 C_s C_u C_{tr}^2 h^{2(r-1)} \sum_{i=1}^N \|v\|^2_{H^r(K_i)} = C_3 h^{2(r-1)} \|v\|^2_{H^r(Q)},$$

where $C_3 = 4 C_s C_u C_{tr}^2$. Using the inequality $|\{a\}^{up}| \leq |a_i| + |a_j|$ and following the proof above, we estimate the next term as follows

$$\sum_{F_{ij} \in \mathcal{F}_I} h \|\{\nabla_x(v - \Pi_h v)\}^{up}\|^2_{L_2(F_{ij})}$$

$$\leq 2 C_u C_{tr}^2 \sum_{i=1}^N h_{K_i} h_{K_i}^{-1} \left( \|\nabla_x(v - \Pi_{h,i} v)\|^2_{L_2(K_i)} + h_{K_i}^2 |\nabla_x(v - \Pi_{h,i} v)|^2_{H^1(K_i)} \right)$$

$$\leq 4C_s C_u C_{tr}^2 h^{2(r-1)} \sum_{i=1}^{N} \|v\|_{H^r(K_i)}^2 = C_4 h^{2(r-1)} \|v\|_{H^r(Q)}^2,$$

where $C_4 = 4C_s C_u C_{tr}^2$. We obtain the final estimate following the techniques in the proof above as follows

$$\sum_{F_{ij} \in \mathcal{F}_I} \delta\theta h \|[\![\partial_t(v - \Pi_h v)]\!]_x\|_{L_2(F_{ij})}^2$$

$$\leq 2\delta\theta C_u C_{tr}^2 \sum_{i=1}^{N} h_{K_i} h_{K_i}^{-1} \left( \|\partial_t(v - \Pi_{h,i} v)\|_{L_2(K_i)}^2 + h_{K_i}^2 |\partial_t(v - \Pi_{h,i} v)|_{H^1(K_i)}^2 \right)$$

$$\leq 4\delta\theta C_u C_{tr}^2 C_s h^{2(r-1)} \sum_{i=1}^{N} \|v\|_{H^r(K_i)}^2 = C_5 h^{2(r-1)} \|v\|_{H^r(Q)}^2,$$

where $C_5 = 4\delta\theta C_u C_{tr}^2 C_s$. □

**Lemma 9** *Let $d$ and $s$ be positive integers with $2 \leq s \leq p+1, s > (d+1)/2$ and let $v \in H^s(Q) \cap H_{0,\underline{0}}^{1,1}(Q)$. Then, there exists a projection $\Pi_h$ from $H_{0,\underline{0}}^{1,1}(Q) \cap H^s(Q)$ to $V_{0h}$ and generic positive constants $C_6$ and $C_7$ such that*

$$\|v - \Pi_h v\|_h \leq C_6 h^{r-1} \|v\|_{H^r(Q)}, \tag{36}$$

$$\|v - \Pi_h v\|_{h,*} \leq C_7 h^{r-1} \|v\|_{H^r(Q)}, \tag{37}$$

*where $h$ is the mesh-size in the physical domain, $r = \min\{s, p+1\}$, $p$ denotes the underlying polynomial degree, and the generic constants $C_6$ and $C_7$ do not depend on $h$ and $v$.*

*Proof* Using the discrete norms as defined in (22) and (26) together with Lemma 8, we complete the proof of the statement. □

Finally, we present the main result for the article, namely the *a priori* error discretization error estimate in the discrete norm $\|\cdot\|_h$.

**Theorem 2** *Let $s$ and $p$ be positive integers with $2 \leq s \leq p+1$ and $s > (d+1)/2$. Further, let $u \in H_{0,\underline{0}}^{1,0}(Q) \cap H^s(Q)$ be the exact solution of our model problem (2), and let $u_h \in V_{0h}$ be the solution to the FEM scheme (13). Then the discretization error estimate*

$$\|u - u_h\|_h \leq C h^{r-1} \|u\|_{H^r(Q)}, \tag{38}$$

*holds, where $C$ is a generic positive constant, $r = \min\{s, p+1\}$, and $p$ denotes the underlying polynomial degree.*

*Proof* By using the coercivity result of Lemma 5, the Galerkin orthogonality (17) and the boundedness of the discrete bilinear form, i.e. Lemma 6, we can derive the following estimates

$$\mu_c \|\Pi_h u - u_h\|_h^2 \leq a_h(\Pi_h u - u_h, \Pi_h u - u_h) = a_h(\Pi_h u - u, \Pi_h u - u_h)$$

$$\leq \mu_b \|\Pi_h u - u\|_{h,*} \|\Pi_h u - u_h\|_h.$$

Hence, we have

$$\|\Pi_h u - u_h\|_h \leq (\mu_b/\mu_c) \|\Pi_h u - u\|_{h,*}. \tag{39}$$

Using (39) together with the estimates (37) and (36) from Lemma 9, we have

$$\begin{aligned}\|u - u_h\|_h &\leq \|u - \Pi_h u\|_h + \|\Pi_h u - u_h\|_h \\ &\leq \|u - \Pi_h u\|_h + (\mu_b/\mu_c)\|\Pi_h u - u\|_{h,*} \\ &\leq (C_6 + C_7(\mu_b/\mu_c)) h^{r-1} \|u\|_{H^r(Q)},\end{aligned}$$

with $C = (C_6 + C_7(\mu_b/\mu_c))$. □

*Remark 4* In practice, the assumption $s > (d+1)/2$ imposed on the regularity is restrictive. In such a case, the domain decomposition of the space-time cylinder is required as presented in [11,12] for the space-time discontinuous Galerkin (dG) using Non uniform rational B-splines (NURBS). Also, the error estimate in Theorem 2 can be further analyzed for the mesh element as follows

$$\|u - u_h\|_h \leq C \left( \sum_{K \in \mathcal{K}_h} h_K^{2(r-1)} \|u\|_{H^r(K)}^2 \right)^{1/2}, \tag{40}$$

where $h_K$ is the mesh size of the mesh element $K$. Such an error estimate (40) is particularly relevant for developing Adaptive Finite Element Method (AFEM) space-time scheme.

## 5 Numerical Results

We present in this section some numerical experiments to validate the theoretical estimates presented in Section 4. The examples considered are motivated by [10, Section 5], where the authors used Isogeometric Analysis (IgA). The numerical results have been performed in FreeFem++ [5]. We will consider a deforming space-time domain with polygonal surface area $\Sigma$ as illustrated in Figure 2. By using unstructured mesh with $p = 1$ and $p = 2$ polynomial degree continuous finite elements, we compute the error in the discrete norm $\|\cdot\|_h$. The resulting linear system is solved using GMRES without preconditioning. We chose the positive parameters $\theta = 0.1$ and $\delta = 10$ in all our numerical experiments.

5.1 Example I

We consider the space-time domain $Q := \Omega(t) \times [0,1] \subset \mathbb{R}^2$ where $\Omega(t) = \{x \in \mathbb{R} : a(t) < x < b(t)\}$ with $t = (0,1), a(t) = -t/2$ and $b(t) = 1 + t/2$, see Figure 2. The exact solution for the model problem (1) is given by $u(x,t) = \sin(\pi x)\sin(\pi t)$. The volume density $f$ and the initial data $u_0 = 0$ on $\Sigma$ are resolved accordingly. By successive mesh refinements using piecewise linear ($p = 1$) and piecewise quadratic ($p = 2$) continuous finite element, the convergence behavior with respect to the discrete norm $\|\cdot\|_h$ can be seen in Figure 3 (right). The solution contours $u_h$ can be seen in Figure 3. Since the exact solution is smooth, we expect optimal convergence rate $\mathcal{O}(h^p)$ for polynomial degrees $p \geq 1$ as predicted by Theorem 2.

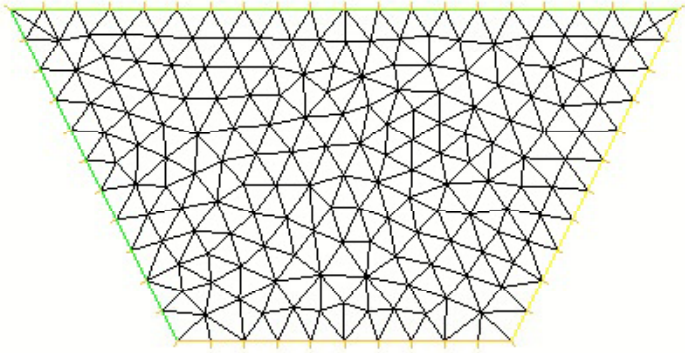

**Fig. 2** Moving spatial domains with underlying unstructured mesh.

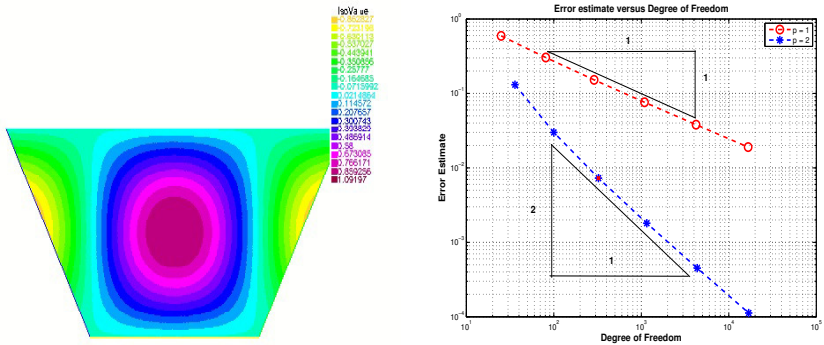

**Fig. 3** Solution contours in the space-time cylinder $Q$ (left) and the plot of the error estimates for degrees $p = 1$ and $p = 2$.

5.2 Example II

We consider again the space-time domain $Q := \Omega(t) \times [0, 1] \subset \mathbb{R}^2$ as described in Example 5.1. In this example, the exact solution for the model problem (1) is given by $u(x, t) = (1 - t)^{1/2} \sin(\pi x)$ and has a singularity at $t = 1$.

The volume density $f$ and the initial data $u_0$ are resolved accordingly. The exact solution $u(x,t)$ belongs to the space $H^{1-\varepsilon}(Q)$ for any $\varepsilon > 0$. By successive mesh refinements using piecewise linear ($p = 1$) and piecewise quadratic ($p = 2$) continuous finite element, the convergence behavior with respect to the discrete norm $\|\cdot\|_h$ can be seen in Figure 4 (right). The solution contours $u_h$ can be seen in Figure 4 (left). Using Theorem 2, we expect a convergence rate of zero i.e. $\mathcal{O}(h^0)$. However, we observe a convergence rate of $\mathcal{O}(h^{0.5})$ for polynomial degrees $p \geq 1$. This is because the solution has a singularity only in time and a full regularity in space. For the linear polynomial degree ($p = 1$), we obtain a convergence rate of 0.75 while for the quadratic polynomial degree ($p = 2$), we obtain the expected convergence rate of 0.5. For the quadratic polynomial degree, we are in the asymptotic range faster than the linear polynomial degree, thus yielding the expected rate in the numerical experiments. The numerical results confirm similar results obtained by using space-time discontinuous Galerkin FEM as presented in [13, Example 2.3.3].

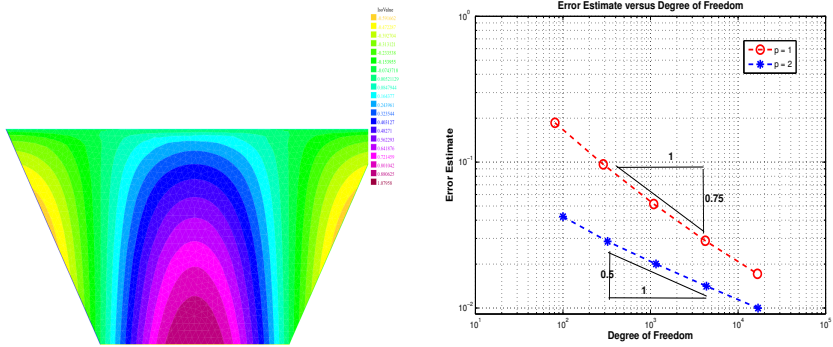

**Fig. 4** Solution contours in the space-time cylinder $Q$ (left) and the plot of the error estimates for degrees $p = 1$ and $p = 2$.

**Conclusion**

In this article, we have presented a new stable space-time finite element scheme for parabolic evolution problems on arbitrary space-time domains. We have presented *a priori* error estimates and numerical examples in the space-time domain $Q \subset \mathbb{R}^2$ using continuous finite element approximations with polynomial degrees $p = 1$ and $p = 2$. For the smooth solution $u$ considered in Subsection 5.1, we have observed the optimal convergence rates as predicted by Theorem 2. However, for non-smooth input data, the solutions have less regularity and thus the *a priori* error estimate obtained in Theorem 2 does not yield the optimal convergence rate or cannot be applied as observed in Subsection 5.2. An appropriate understanding of the behavior of such non-smooth solutions will be to consider anisotropic Sobolev spaces. The new results present the possibility for space-time adaptivity in moving domains, also the possibility to use preconditioning, fast solvers like multi-grid and domain decomposition solvers.

**Acknowledgment**